\documentclass[12pt]{article}
 
\usepackage{coursE}

\title{About the $L^2$ analyticity of Markov operators on graphs}

\author{{\Large Joseph {\sc Feneuil}} \\
School of Mathematics, University of Minnesota, \\
 127 Vincent Hall, 206 Church St. SE, \\
Minneapolis, Minnesota 55455 \\
jfeneuil@umn.edu \\
\vspace{1cm}}
\date{\today}

\DeclareMathOperator{\Sp}{Sp}

\begin{document}

\maketitle

\begin{abstract}
Let $\Gamma$ be a graph and $P$ be a reversible random walk on $\Gamma$. 
From the $L^2$ analyticity of the Markov operator $P$, we deduce that an iterate of odd exponent of $P$ is `lazy', that is there exists an integer $k$ such that the transition probability (for the random walk $P^{2k+1}$) from a vertex $x$ to itself is uniformly bounded from below.
The proof does not require the doubling property on $\Gamma$ but only a polynomial control of the volume.
\end{abstract}

\tableofcontents

\pagebreak

\section{Introduction and statement of the results}

In 1990, in \cite[Section III]{CSC}, Coulhon and Saloff-Coste introduced the following notion of analyticity for discrete semigroups: 
let $(X,m)$ be a space equipped with a measure. We say that an operator $T$ on $L^p(X,m)$ - or the discrete semigroup $(T^k)_{k\in \N}$ - is {\bf analytic} if there exists $C>0$ such that for any $k\in \N^*$,
\begin{equation} \label{DiscreteAnalyticity}
\|(I-T)T^k\|_{p\to p} = \|T^k-T^{k+1}\|_{p\to p} \leq \frac{C}{k},
\end{equation}
where $\|L\|_{p\to p} = \sup_{\|f\|_p = 1} \|L\, f\|_p$. 
One of the characterizations of the $L^p$-analyticity for a (continuous) semigroup $H_t$ with generator $A$ is the existence of $C>0$ such that for any $t>0$,
\begin{equation}  \label{ContinuousAnalyticity}
\|AH_t\|_{p\to p} = \|\dr_t H_t \|_{p\to p} \leq \frac{C}{t}.
\end{equation}
If the discrete derivative of a map $f$ defined on $\N$ is $\dr_k f(k) =f(k+1)-f(k)$, Condition \eqref{DiscreteAnalyticity} can be thus seen as the discrete analogue of Condition \eqref{ContinuousAnalyticity}.

\medskip

Comparison with the continuous notion of analyticity gives that $-\Delta = T-I$ should play the role of the `generator' of the discrete semigroup. So if $T$ is, say, analytic and self-adjoint, then $\Delta = I-T$ can be chosen as the positive Laplacian and can be used to define fractional Sobolev spaces or Besov spaces with appropriate embedding inequalities (see \cite[Section III]{CSC}, with methods of \cite{CSC2}). Additional links of discrete analyticity with Littlewood-Paley theory can also be found in \cite{Dungey2,LMX}. Discrete analyticity is also used to establish time derivative estimates on the discrete time heat kernel (see \cite{Blunck00,Dungey}, see \cite{christ} for the first but more complicated proof of this fact): if $(T^k)_{k\in \N}$ is analytic, power bounded on $L^2$ and if there exist $C,c>0$ such that the kernel of $T^k$ satisfies
$$h_{k}(x,y) \leq Ck^{-d/\beta} \exp\left(-c\left(\frac{d(x,y)^\beta}{k}\right)^{\beta/(\beta-1)}\right) \qquad \forall x,y\in X, \, \forall k\in \N^*$$ 
for some $d>0$ and $\beta \geq 2$, then there exist $C,c>0$ such that the kernel of $T^k$ satisfies
$$|h_{k+1}(x,y) - h_k(x,y)| \leq Ck^{-\frac{d}{\beta} - 1} \exp\left(-c\left(\frac{d(x,y)^\beta}{k}\right)^{\beta/(\beta-1)}\right) \qquad \forall x,y\in X, \, \forall k\in \N^*.$$
As another application, let us cite that discrete analyticity is a necessary condition for the maximal regularity of the discrete time evolution equation, which satisfies $u_{k+1} + Tu_k = f_k$ for all $k\in \N$ and the initial condition $u_0=0$ (see \cite{Blunck01}).

\medskip

In the sequel, the space we consider is a graph $\Gamma$ (see Subsection \ref{defgraphs} for complete definitions) and $P$ denotes a reversible random walk (in particular, $P$ is self-adjoint and $P$ is a contraction on $L^1$ and $L^\infty$). These conditions are not enough to insure the discrete analyticity of $P$. Indeed, if $P$ is the standard random walk on $\Z$, that is the transition probability from $x$ to $y$ is $\frac12$ if $|x-y|=1$ and 0 otherwise, then $P$ is not analytic on $L^2$. 

We want characterizations of \eqref{DiscreteAnalyticity}. Denote by $\Sp_{L^q}(P)$ the spectrum of $P$ in $L^q(\Gamma)$ and by $\mathbb D$ the open unit disk in $\C$. A first characterization uses a result of Nevanlinna that runs in a very general context: $P$ is analytic (and power bounded) on $L^q(\Gamma)$ if and only if $\Sp_{L^q}(P) \subset \mathbb D \cup \{1\}$ and the semigroup $(e^{-t(I-P)})_{t>0}$ is $L^q$ analytic (see \cite[Theorem 4.5.4]{Nevanlinna1}, \cite[Theorem 2.1]{Nevanlinna2} or \cite{Blunck00}, see also \cite{LeMerdy} and the references therein for the spectral properties of Ritt operators).

If $P$ is selfadjoint and $P$ is a contraction on $L^2$, we can obtain another criterion of discrete analyticity via the spectral theorem. Indeed, under these hypotheses,
$P$ is a contraction on $L^2$, $\Sp_{L^2}(P) \subset [-1,1]$ and thus $P$ is $L^2$ analytic if and only if $-1$ doesn't belong to $\Sp_{L^2}(P)$. 
Then, for $1<p<+\infty$, the $L^p$ analyticity of $P$ can be deduced from its $L^2$ analyticity (as explained at the end of \cite[Section III]{CSC}, by using an interpolation argument by Stein found in \cite{Stein1}).

Although $-1 \notin \Sp_{L^2}(P)$ is equivalent to the $L^2$ analyticity, another stronger condition is found in the literature:
we say that the operator $P$ (or the semigroup $P^k$) satisfies \eqref{LB} if the transition probability from one point $x$ to itself is uniformly bounded from below by a positive constant (see subsection \ref{ss1.2} for mathematical statement). Such a random walk is sometimes called `lazy'. 
Condition \eqref{LB} can be seen in several articles \cite{HebSC,Russ, CGZ,Dungey,Fen1,CCFR}, and the fact that \eqref{LB} implies \eqref{DiscreteAnalyticity} can be found for instance in \cite[Lemma 1.3]{Dungey}. 
 
\medskip

The purpose of this paper is to provide a converse to the fact that \eqref{LB} implies \eqref{DiscreteAnalyticity}. More precisely, when the graph $\Gamma$ has at most polynomial growth, we proved that the discrete analyticity of $P$ forces an iterate of odd exponent of $P$ to satisfy \eqref{LB}. The last characterization of \eqref{DiscreteAnalyticity} given in Theorem \ref{MainTheo} is a practical one, which can be used to check quickly the analyticity of a given graph. As an application, we slightly improve a result in \cite{Russ,Fen4,CCFR} that establishes the $L^p$ boundedness of the Riesz transform for $1<p<2$ on graphs under Gaussian or sub-Gaussian estimates.

\medskip

The remainder of the first chapter is devoted to, in order, the presentation of the graph structure, the introduction of the assumptions used in this article, 
and the accurate statement of our main result. The second section is dedicated to the proof of the main result. 
Finally, in the last section, we present an application of our main result to weaken the assumptions used in \cite{Russ,CCFR}. 

\medskip

\noindent {\bf Acknowledgement:} The author was supported by the ANR project ``Harmonic Analysis at its Boundaries'',   ANR-12-BS01-0013-03. The author would also like to thank the referee for pointing out some references given here and for comments that make the article easier to read.

\subsection{The setting} \label{defgraphs}

Let $\Gamma$ be a countable (infinite) set and $\mu_{xy} = \mu_{yx} \geq 0$ a symmetric weight on $\Gamma \times \Gamma$. The set $E$ of edges is $\{(x,y)\in \Gamma^2, \, \mu_{xy}>0\}$. We say that $x$ and $y$ are neighbors - and we write $x\sim y$ - if $(x,y)\in E$. 
The graph is outfitted with the classical distance. A sequence of vertices $x_0,\dots,x_n$ is a path (of length $n$) if for all $i\in \{1,n\}$, $x_{i-1} \sim x_i$. The distance $d$ between $x$ and $y$ is then the length of the shortest path linking $x$ to $y$. We assume that the set is connected, that is $d(x,y)<+\infty$ for any $x,y\in \Gamma$. Moreover, we write $B(x,r)$ for the set $\{y\in \Gamma, \, d(x,y)<r\}$.

We define the weight $m(x)$ of a vertex $x \in \Gamma$ by $m(x) = \sum_{x\sim y} \mu_{xy}$. Note that since $\Gamma$ is connected, $m(x)>0$ for any $x\in \Gamma$.
More generally, the volume (or measure) of a subset $E \subset \Gamma$ is defined as $m(E) := \sum_{x\in E} m(x)$ and the term $V(x,r)$ is used for $m(B(x,r))$.$L^p(\Gamma)$ spaces are defined as follows: for all $1\leq p < +\infty$, a function $f$ on $\Gamma$ belongs to $L^p(\Gamma,m)$ (or $L^p(\Gamma)$) if
$$\|f\|_p := \left( \sum_{x\in \Gamma} |f(x)|^p m(x) \right)^{\frac{1}{p}} < +\infty.$$

The weight $\mu_{xy}$ on the edges is also used to define a random walk. For all $x,y\in \Gamma$ the discrete-time reversible Markov kernel $p$ associated with the measure $m$ by $p(x,y) = \frac{\mu_{xy}}{m(x)m(y)}$.
The discrete kernel $p_k(x,y)$ is then defined recursively for all $k\geq 0$ by
\begin{equation}
\left\{ 
\begin{array}{l}
p_0(x,y) = \frac{\delta(x,y)}{m(y)} \\
p_{k+1}(x,y) = \sum_{z\in \Gamma} p(x,z)p_k(z,y)m(z).
\end{array}
\right.
\end{equation}
Notice that for all $k\geq 1$, 
the kernel $p_k$ is symmetric, that is $p_k(x,y) = p_k(y,x)$ for all $x,y\in \Gamma$. Moreover, $p_k(x,y)m(y)$ denotes the probability to go from $x$ to $y$ in $k$ steps.
For all functions $f$ on $\Gamma$, we define $P$ as the operator with kernel $p$, that is, for any $x\in \Gamma$ and any function $f$,
$Pf(x) = \sum_{y\in \Gamma} p(x,y)f(y)m(y).$
It is easily checked that $p_k$ is the kernel of the operator $P^k$ with respect to the measure $m$.Since $p_k(x,y) \geq 0$ and $\sum_{y\in \Gamma} p_k(x,y)m(y) = 1$, $P^k$ contracts $L^p(\Gamma)$ for all $k\in \N$ and all $p\in [1,+\infty]$.

\medskip

We need for the sequel to define additional metrics on the graph. Let $\epsilon>0$. 
We say that $x,y \in \Gamma$ are $\epsilon$-neighbors if $p(x,y)\min\{m(y),m(x)\} >\epsilon$. If $x$ and $y$ are $\epsilon$-neighbors, we write $x\sim_\epsilon y$, otherwise, the notation $x\not\sim_\epsilon y$ is used. An $\epsilon$-path (joining $x$ to $y$, of length $n$) in $\Gamma$ is a sequence $x=x_0,x_1,\dots,x_n = y$ of vertices such that $x_{i-1}\sim_\epsilon x_i$ for all $i\in \bb 1,n\bn$.  
A closed $\epsilon$-path is a path $x=x_0,x_1,\dots,x_n =x$ joining $x$ to itself.
The distance $d_\epsilon(x,y)$ between $x$ and $y$ is defined as the length of the shortest $\epsilon$-path linking $x$ to $y$. 
We write $d_\epsilon(x,y)=+\infty$ if no $\epsilon$-path linking $x$ to $y$ exists.
The $\epsilon$-ball $B_\epsilon(x,r)$ denotes $\{y\in \Gamma, \, d_\epsilon(x,y)<r\}$, and $V_\epsilon(x,r)$ is the quantity $m(B_\epsilon(x,r))$. Remark in particular that $B_\epsilon(x,r) \subset B(x,r)$ and thus $V_\epsilon(x,r) \leq V(x,r)$.

Finally, let us state the definition of the boundary of a set: the $\epsilon$-boundary of $E\subset \Gamma$ is the set 
$$\partial_\epsilon E = \{x\in E, d_\epsilon(x,E^c) = 1\}.$$

We assume in all the article that $\Gamma$ is locally uniformly finite, that is there exists $N_{dvl} \in \N$ such that, for all $x\in \Gamma$,
$$\#B_0(x,2) =\#\left(\{y\in \Gamma, \, y\sim x\} \cup \{x\} \right) \leq N_{dvl}.$$

\subsection{Statement of the assumptions and main result}

\label{ss1.2}

\begin{defi}
 Let $k\in \N^*$. We say that $P^k$ satisfies \eqref{LB} if there exists $\epsilon_{LB} >0$ (that can depend on $k$) such that
\begin{equation} \label{LB} \tag{LB}
p_k(x,x)m(x) \geq \epsilon_{LB} \qquad \forall x\in \Gamma.
\end{equation}
\end{defi}

\begin{defi}
We say that $(\Gamma,\mu)$ satisfies \eqref{PDV} if there exists $C_{dvp},D >0$ such that
\begin{equation} \label{PDV} \tag{DVP}
V(x, r) \leq C_{dvp} r^D m(x) \qquad \forall x\in \Gamma, \, \forall r\in \N^*.
\end{equation}
\end{defi}

\begin{rmk}
Observe that the property \eqref{PDV} is weaker than the doubling property \eqref{DV}, that states that there exists $C_{dv}>0$ such that
\begin{equation} \label{DV} \tag{DV}
V(x, 2 r) \leq C_{dv} V(x,r) \qquad \forall x\in \Gamma, \, \forall r\in \N^*.
\end{equation}
Furthermore, the property \eqref{PDV} is stronger than the local doubling property \eqref{DVL}, which says that there exists $C_{dvl}>0$ such that
\begin{equation} \label{DVL} \tag{DVL}
V(x,2) \leq C_{dvl} m(x) \qquad \forall x\in \Gamma.
\end{equation}
\end{rmk}

\begin{rmk} \label{P2LB}
Observe that if $(\Gamma,\mu)$ satisfies \eqref{DVL}, then, $P^2$ satisfies \eqref{LB} (see for example \cite[Lemma 3.2]{CGZ}).
\end{rmk}

Here and after, if $E\subset \Gamma$, by ``partition'' of $E$, we mean a collection of pairwise (possibly empty) disjoint subsets of $E$.
Let us now state the following result, proven in subsection \ref{par2pathed}:

\begin{prop} \label{2pathedprop}
Let $\epsilon > 0$ and $E \subset \Gamma$. The following conditions are equivalent
\begin{enumerate}[(i)]
 \item all closed $\epsilon$-paths in $E$ are of even length,
 \item there exists a partition $(E_0,E_1)$ of $E$ such that if $i\in \{0,1\}$ then 
\begin{equation} \label{iiProp2pathed} x,y\in E_i \f x\not\sim_\epsilon y.\end{equation}
\end{enumerate}
A set satisfying one of the two equivalent conditions is said to be $\epsilon$-bipartite.
\end{prop}

\begin{rmk}
If $E$ is $\epsilon$-bipartite, then for all $x\in E$, $p(x,x)m(x) \leq \epsilon$.
\end{rmk}

Here is our main result:

\begin{theo} \label{MainTheo}
Let $(\Gamma,\mu)$ be a graph satisfying \eqref{PDV}. The following conditions are equivalent:
\begin{enumerate}[(i)]
 \item $P$ is analytic on $L^2$, that is there exists $C>0$ such that for any $k\in \N^*$,
  $$\|(I-P)P^k\|_{2\to 2} \leq \frac {C}k,$$
 \item for any $q\in (1,+\infty)$, $P$ is analytic on $L^q$, that is there exists $C_q>0$ such that for any $k\in \N^*$,
 $$\|(I-P)P^k\|_{q\to q} \leq \frac {C_q}k,$$
 \item for all $q\in [1,+\infty]$, $-1 \notin Sp_{L^q}(P)$,
 \item there exists $q\in [1,+\infty)$ such that  $-1 \notin Sp_{L^q}(P)$,
 \item there exists $k\in \N$ such that $P^{2k+1}$ satisfies \eqref{LB}, that is there exist $k\in \N$ and $\epsilon>0$ such that for all $x\in \Gamma$, one has 
 $$p_{2k+1}(x,x) m(x) \geq \epsilon,$$
 \item there exist $\epsilon>0$ and $r_0\in\N^*$ such that if the $\epsilon$-ball $B_\epsilon:=B_\epsilon(x,r)$ is $\epsilon$-bipartite, then $r\leq r_0$.
\end{enumerate}
\end{theo}

\begin{rmk}
Assume that the weights $\mu_{xy}$ on the edges are `locally doubling', that is there exists $C>0$ such that for any vertices $x,y,z\in \Gamma$ satisfying $x\sim y$ and $x\sim z$, we have $\mu_{xy} \leq C \mu_{xz}$. In this case, note that there exists $\epsilon>0$ such that for any $x\sim y$, $p(x,y)m(y) > \epsilon$. Therefore, any path is an $\epsilon$-path and any ball is an $\epsilon$-ball. 

With this additional assumption, item (vi) in Theorem \ref{MainTheo} means that there exists $r_0>0$ such that for any $x\in \Gamma$, we can find a closed path going through $x$ whose length is odd and smaller than $r_0$. It is easy to see that this condition is not satisfied for the standard random walk on $\Z^n$.
\end{rmk}

\section{Proofs}

\subsection{Proposition \ref{2pathedprop}}

\label{par2pathed}

\begin{proof}
Let us prove $(i) \f (ii)$. We assume first that $E$ is $\epsilon$-connected, 
that is, for all $x,y\in E$, $d_\epsilon(x,y)<+\infty$. 
Fix $e\in E$. Define
$$E_0 = \{x\in E, \, \text{there exists an $\epsilon$-path of even length in $E$ linking $e$ to $x$}\}$$
and
$$E_1 = \{x\in E, \, \text{there exists an $\epsilon$-path of odd length in $E$ linking $e$ to $x$}\}.$$
\begin{itemize}
 \item Since $E$ is $\epsilon$-connected, one has $E_0 \cup E_1 = E$.
 \item Let us check that $E_0 \cap E_1 = \emptyset$. Assume the contrary and let $f\in E_0 \cap E_1$. 
Then, since $f\in E_0$, there exist an integer $n$ and an $\epsilon$-path in $E$ given by: 
$$e=x_0,\dots,x_{2n}=f,$$
and since $f\in E_1$, there exist an integer $m$ and an $\epsilon$-path in $E$ given by:
$$e=y_0,\dots,y_{2m+1}=f.$$
Thus $x_0,\dots,x_{2n},y_{2m},\dots,y_0$ is a closed $\epsilon$-path in $E$ of odd length, which is a contradiction.

 \item Let $x,y\in E$ be such that $x\sim_\epsilon y$. If $x\in E_i$, then there exist an integer $n$ and an $\epsilon$-path $e=x_0,\dots,x_{2n+i}=x$.
Consequently, $e=x_0,\dots,x_{2n+i},y$ is an $\epsilon$-path of length $2n+1+i$, and thus $y\in E_{1-i}$.
\end{itemize}

Assume now that $E$ is not $\epsilon$-connected. 
Let $(E^1,\dots,E^n)$ be the $\epsilon$-connected components of $E$ 
(that is, $(E^1,\dots,E^n)$ is a partition of $E$ such that for all $x\in E^j$ and $y\in E^k$, $d(x,y)<+\infty$ if $j=k$ and $d_\epsilon(x,y) = +\infty$ otherwise).
According to the above proof, for all $k\in \{1,\dots,n\}$, there exists a partition $(E_0^k,E_1^k)$ of $E^n$ satisfying \eqref{iiProp2pathed}.
Since no $\epsilon$-edges are linking $E^j$ to $E^k$ whenever $j\neq k$, the partition $(E_0,E_1)$ of $E$, where for any $i\in \{0,1\}$
$$E_i := \bigcup_{k=1}^n E_i^k$$
satisfies \eqref{iiProp2pathed}.

\medskip

Let us now turn to the proof of $(ii) \f (i)$. There exists a partition $(E_0,E_1)$ of $E$ satisfying \eqref{iiProp2pathed}.
Let $e=x_0,\dots,x_n=e$ be a closed $\epsilon$-path in $E$. Without loss of generality, we can assume that $e\in E_0$.
Then, for all $k\in \bb 0,n\bn$, $x_k \in E_i$ if and only if $k-i$ is even. 
As a consequence, since $x_n = e \in E_0$, the integer $n$ is even.
As a conclusion, any closed $\epsilon$-path in $E$ is of even length.
\end{proof}

\subsection{Theorem \ref{MainTheo}}

Here is the plan of the proof. First, we talk about the implications 
$$ (v) \f (iii) \f (ii) \f (i) \f (iv)$$
that can already be found in the literature or that are easy to prove. Then, we prove the implication $(vi) \f (v)$ and we finish by the implication $ (iv) \f (vi)$, which is actually the only one that requires Assumption \eqref{PDV}.

\begin{proof} [Proof of Theorem \ref{MainTheo}, $ (v) \f (iii) \f (ii) \f (i) \f (iv)$]

Let us begin by $(v)$ implies $(iii)$. Let $q\in [1,+\infty]$ and let $k\in \N$ be such that $P^{2k+1}$ satisfies \eqref{LB}. According to \cite[Lemma 1.3]{Dungey},
\begin{equation} \label{Dungeyresult}
\Sp_{L^q}(P^{2k+1}) \subset \mathbb D \cup \{1\} := \{z\in \C, \, |z|< 1\} \cup \{1\}.
\end{equation}
Assume now that $-1\in \Sp_{L^q}(P)$, then $(-1)^{2k+1} = -1 \in \Sp_{L^q}(P^{2k+1})$, which contradicts \eqref{Dungeyresult}.

\smallskip

The implication $(iii) \f (ii)$ is proved in \cite[Section III]{CSC} and the one $(ii) \f (i)$ is immediate.

\smallskip

Assume $(i)$ and we want to prove $(iv)$. By contraposition, it suffices to prove that $-1 \in \Sp_{L^2}(P)$ yields that $P$ is not $L^2$ analytic, or that for any $k\in \N^*$, there exists $f\in L^2(\Gamma)$ such that $\|(I-P)P^k f\|_2 \geq 1$. One has
\[\begin{split}
(I-P)P^k & = [2I - (I+P)] [(P+I)-I]^k = 2 \sum_{j=0}^k \begin{pmatrix} k \\ j \end{pmatrix}  (-1)^{k-j} (P+I)^j - \sum_{j=0}^k \begin{pmatrix} k \\ j \end{pmatrix} (-1)^{k-j} (P+I)^{j+1}  \\
& = 2(-1)^kI - (P+I)  \sum_{j=0}^k \left[\begin{pmatrix} k \\ j \end{pmatrix} +2 \begin{pmatrix} k \\ j+1 \end{pmatrix} \right]   (-1)^{k-j} (P+I)^{j}.
\end{split}\]
The operator 
$$\sum_{j=0}^k \left[\begin{pmatrix} k \\ j \end{pmatrix} +2 \begin{pmatrix} k \\ j+1 \end{pmatrix} \right]   (-1)^{k-j} (P+I)^{j}$$
is bounded on $L^2$ by some constant $C_k$ depending on $k$. As a consequence, for any $f\in L^2(\Gamma)$,
$$\|(I-P)P^k f \|_{2} \geq 2 \|f\|_2 - C_k \|(P+I)f\|_2.$$
Yet, since $-1\in \Sp_{L^2}(P)$, $0$ belongs to $\Sp_{L^2}(P+I)$ and there exists $f_0\in L^2(\Gamma)$ such that $\|f_0\|_{2}=1$ and $\|(P+I)f_0\|_{2} \leq C_k^{-1}$. It follows that
$$\|(I-P)P^k f_0 \|_{2} \geq  \|f_0\|_2 = 1$$
and then that $\|(I-P)P^k \|_{2\to 2} \geq 1$.
\end{proof}

\begin{proof} [Proof of Theorem \ref{MainTheo}, $(vi)$ implies $(v)$]

We prove again the implication by contraposition, that is we assume that for all $k\in \N$ and all $\eta>0$, there exists $x_0 \in \Gamma$ such that 
\begin{equation} \label{x0def} p_{2k+1}(x_0,x_0)m(x_0) \leq \eta .\end{equation} Fix $\epsilon>0$ and $r\in \N^*$. 

\medskip

\noindent {\bf First step:}
We claim that there exists $x_0\in \Gamma$ such that there doesn't exist any $\epsilon$-path of length $2r-1$ starting from $x_0$.\par

Let $x_0$ satisfy \eqref{x0def} for $k=r-1$ and $\eta = \epsilon^{2r-1}$. Choose $x_0,x_1,\dots,x_{2r-1} = x_0$ any path of length $2r-1$. One has
$$\prod_{i\in \bb 1,2r-1\bn} p(x_{i-1},x_{i})m(x_i) \leq p_{2r-1}(x_0,x_0)m(x_0) \leq \epsilon^{2r-1}.$$
Then there exists $i \in \bb 1,2r-1\bn$ such that $p(x_{i-1},x_{i})m(x_i) \leq (\epsilon^{2r-1})^{\frac{1}{2r-1}} = \epsilon$ and the path is not an $\epsilon$-path.

\medskip

\noindent  {\bf Second step:}
Define now
$$E_0 = \{x\in B_\epsilon(x_0,r), \, d_{\epsilon}(x_0,x) \text{ is even}\}$$
and
$$E_1 = \{x\in B_\epsilon(x_0,r), \, d_{\epsilon}(x_0,x) \text{ is odd}\}.$$

The couple $(E_0,E_1)$ forms a partition of $B_\epsilon(x_0,r)$. Let $x,y \in B_\epsilon (x_0,r)$ be such that $x\sim_\epsilon y$. 
Then $|d_\epsilon(x,x_0)-d_\epsilon(x_0,y)|\leq 1$. 
We claim that $|d_\epsilon(x,x_0)-d_\epsilon(x_0,y)|=1$. Assume on the contrary that $d_\epsilon(x_0,x)=d_\epsilon(x_0,y) = n <r$. 
Then there exist two $\epsilon$-paths $x_0,\dots,x_n=x$ and $x_0,y_1,\dots,y_n=y$.
Thus $x_0,\dots,x_n,y_n,\dots,y_1,x_0$ is an $\epsilon$-path of length $2n+1$. 
We want to extend this closed $\epsilon$-path to a closed $\epsilon$-path of length $2r-1$.
Two cases may happen. 
Either $n=0$ and then, since $x_0=x=y$, $x_0,x_0$ is an $\epsilon$-path of length $1$. Hence
$$\underbrace{x_0,x_0,\dots,x_0}_{2r \text{ times}}$$
is an $\epsilon$-path of length $2r-1$. Or $n>0$ and $x_0,x_1,x_0$ is an $\epsilon$-path of length $2$ and thus
$$\underbrace{x_0,x_1,\dots,x_0,x_1}_{r-n-1 \text{ times }},x_0,\dots,x_n,y_n,\dots,y_1,x_0$$
is an $\epsilon$-path of length $2r-1$. By contradiction with the first step, one has then the desired conclusion:
\begin{equation} \label{differentparity}|d_\epsilon(x,x_0)-d_\epsilon(x_0,y)|= 1.\end{equation}
The identity \eqref{differentparity} yields that $d_\epsilon(x,x_0)$ and $d_\epsilon(x_0,y)$ don't have the same parity, 
and thus that $x$ and $y$ do not belong to the same $E_i$.

\medskip

\noindent  {\bf Conclusion:}
A way to rephrase the second step is: $(E_0,E_1)$ is a partition of $B(x_0,r)$ such that if $x,y\in E_i$ for some $i\in \{0,1\}$, then $x\not\sim_\epsilon y$. Proposition \ref{2pathedprop} gives then that $B(x_0,r)$ is $\epsilon$-bipartite.
\end{proof}

We turn to the last implication of the proof of Theorem \ref{MainTheo}. Let us start with a definition.

\begin{defi}
Let $\epsilon\geq 0$. Say that a collection of $\epsilon$-balls $(B_\nu)_{\nu\in \N^*}$ has the property \eqref{NB} if and only if, for all $\nu\geq 1$,
\begin{equation} \label{NB} \tag{NB} 
 m(\dr_\epsilon B_\nu) \leq \frac{1}{\nu} m(B_\nu).
\end{equation}
\end{defi}

\begin{prop} \label{PropNB}
Let $(\Gamma,\mu)$ be a weighted graph satisfying \eqref{PDV} and $\epsilon \geq 0$. Let $(\bar B_n)_{n\in \N^*}$ be a collection of $\epsilon$-balls such that, for all $n\geq 1$, the radius of $\bar B_n$ is equal to $n$. 
Then there exists a collection $(B_\nu)_{\nu\in \N^*}$ of $\epsilon$-balls that satisfies:
\begin{enumerate}[(i)]
 \item for all $\nu\in \N^*$, there exists $n \in \N^*$ such that $B_\nu \subset \bar B_n$,
 \item $(B_\nu)_{\nu\in \N^*}$ has the property \eqref{NB}.
\end{enumerate}
\end{prop}

\begin{proof}
Define $x_n$ as the center of the $\epsilon$-ball $\bar B_{n}$. Assume that there exists $K>0$ such that for all $n\in \N^*$ and all $l\in \bb 1,n-1\bn$,
\begin{equation} \label{Contradiction} m(\dr_\epsilon B_\epsilon(x_n,l)) \geq K m(B_\epsilon(x_n,l))\end{equation}
Since we have the inclusion
 $$B_\epsilon(x_n,l+1) \supset \dr_\epsilon B_\epsilon(x_n,l+1) \cup B_\epsilon(x_n,l),$$
one has for all $n\in \N$ and all $l\in \bb 1,n-1\bn$,
$$ m(B_\epsilon(x_n,l+1)) \geq (1+K) m(B_\epsilon(x_n,l)),$$
and thus by induction on $l$, one obtains for all $n\in \N^*$
$$ m(\bar B_n) = m(B_\epsilon(x_n,n)) \geq (1+K)^{n-1} m(x_n).$$
Yet, \eqref{PDV} yields that there exist $C,D>0$ such that, for all $n\in\N^*$,
$$ m(\bar B_{n}) \leq C n^{D} m(x_n).$$
Hence, we obtain for all $n\in \N^*$
$$(1+K)^{n-1} \leq C n^D,$$
which is impossible. Therefore, \eqref{Contradiction} is false, that is, for all $\nu\in \N^*$, there exist $n_\nu \in \N^*$, $l_\nu\in \bb 1,n_\nu-1\bn$ such that
$$m(\dr_\epsilon B_\epsilon(x_{n_\nu},l_\nu)) \leq \frac{1}{\nu} m(B_\epsilon(x_{n_\nu},l_\nu)),$$
and the collection $(B_\nu)_{\nu\in \N^*}$ where $B_\nu = B_\epsilon(x_{n_\nu},l_\nu) \subset \bar B_{n_\nu}$ satisfies (i) and (ii).
\end{proof}

We are now ready to prove the last part of Theorem \ref{MainTheo}.

\begin{proof} [Proof of Theorem \ref{MainTheo}, (iv) implies  (vi)]

We prove as before the result by contraposition. Assume $\neg (v)$, that is, for any $r\in \N^*$ and any $\epsilon>0$ there exists $x_0\in \Gamma$ such that $B_\epsilon(x_0,r)$ is $\epsilon$-bipartite.

Fix $\eta >0$ and $q\in [1,+\infty)$. In order to prove $\neg (vi)$, it suffices to show the existence of a nonzero function $f_\eta$ on $\Gamma$ such that 
\begin{equation} \label{I+Pfeta}
\|(I+P)f_\eta\|_{q}^q \leq \eta \|f_\eta\|_{q}^q.
\end{equation} 
Let $\epsilon> 0$ to be fixed later. By assumption, there exists a collection of $\epsilon$-balls $(\bar B_r^\epsilon)_{r\in \N^*}$ 
such that for all $r\in \N^*$, $\bar B_r^\epsilon$ is of radius $r$ and $\bar B_r^\epsilon$ is $\epsilon$-bipartite.

\noindent According to Proposition \ref{PropNB}, there exists then a collection $(B_\nu^\epsilon)_{\nu\in\N^*}$ of $\epsilon$-balls that satisfies
\begin{enumerate}
 \item  $\ds m(\dr_\epsilon B_\nu^\epsilon) \leq \frac{1}{\nu} m(B_\nu^\epsilon)$,
 \item  for all $\nu\in \N^*$, there exists $r_\nu \in \N^*$ such that $B_\nu^\epsilon \subset \bar B_{r_\nu}$.  
Hence, for all $\nu\in \N^*$, $B_\nu^\epsilon$ is $\epsilon$-bipartite.
\end{enumerate}

Let $\nu \in \N^\ast$ to be also fixed later. Since $B_\nu^\epsilon$ is $\epsilon$-bipartite, Proposition \ref{2pathedprop} provides a partition ($E_0$, $E_1$) of $B_\nu^\epsilon$. 
Define $f_\eta:\Gamma\rightarrow \R$ by
$$f_\eta(x) = \left\{\begin{array}{ll} 0 & \text{ if } x\notin B_\nu^\epsilon \\ 1 & \text{ if } x\in E_0 \\ -1 & \text{ if } x\in E_1.
\end{array}\right. $$
Then one has
\begin{equation} \label{1star} \begin{split} 
 \|(I+P)f_\eta\|_{q}^q  & = \sum_{x\in \Gamma} \left| \sum_{y\in \Gamma}[f_\eta(x)+f_\eta(y)]  p(x,y) m(y) \right|^q m(x) \\
&  \leq \sum_{x\in \Gamma}  \sum_{y\in \Gamma}|f_\eta(x)+f_\eta(y)|^q  p(x,y) m(y)  m(x) \\
&  \leq \sum_{\begin{subarray}{c} x,y\in \Gamma \\ y \sim_\epsilon x \end{subarray}}|f_\eta(x)+f_\eta(y)|^q  p(x,y) m(y) m(x) + \sum_{\begin{subarray}{c} x,y\in \Gamma \\ y \not\sim_\epsilon x \end{subarray}}|f_\eta(x)+f_\eta(y)|^q  p(x,y) m(y) m(x).
  \end{split} \end{equation}
Let us estimate the first term. Remark that $|f_\eta(x)+f_\eta(y)|=0$ when either $x\sim_\epsilon y$ and $x,y \in B^\epsilon_\nu$ or $x,y \in (B^\epsilon_\nu)^c$. Indeed, if  $x,y \in (B^\epsilon_\nu)^c$, then the result is obvious and if $x\sim_\epsilon y$ and $x,y \in B^\epsilon_\nu$, then either $(x,y)\in E_0\times E_1$ or $(x,y) \in E_1 \times E_0$ and $f_\eta(x)+f_\eta(y) = 1-1 = 0$. Moreover, note that $|f_\eta(x)+f_\eta(y)|=1$ when $(x,y)\in (B_\nu^\epsilon)^c \times B_\nu^\epsilon$. Therefore, we have
\begin{equation} \label{2star} \begin{split} 
& \sum_{\begin{subarray}{c} x,y\in \Gamma \\ y \sim_\epsilon x \end{subarray}} |f_\eta(x)+f_\eta(y)|^q  p(x,y) m(y) m(x) \\
& \leq \sum_{y \in \dr_\epsilon B_\nu^\epsilon} \sum_{x\in (B_\nu^\epsilon)^c} |f_\eta(x)+f_\eta(y)|^q  p(x,y) m(y) m(x) +  \sum_{x \in \dr_\epsilon B_\nu^\epsilon} \sum_{y\in (B_\nu^\epsilon)^c} |f_\eta(x)+f_\eta(y)|^q  p(x,y) m(y) m(x) \\
& \quad = 2 \sum_{y \in \dr_\epsilon B_\nu^\epsilon} \sum_{x\in (B_\nu^\epsilon)^c} |f_\eta(x)+f_\eta(y)|^q  p(x,y) m(y) m(x) \\
& \leq  2 \sum_{y \in \dr_\epsilon B_\nu^\epsilon} \sum_{x\in (B_\nu^\epsilon)^c}  p(x,y) m(y) m(x) \leq 2 \sum_{y \in \dr_\epsilon B_\nu^\epsilon} \left(\sum_{x\in \Gamma}  p(x,y) m(x) \right) m(y) = 2 m(\dr_\epsilon B^\epsilon_\nu) \\
& \leq \frac{2}{\nu} m(B^\epsilon_\nu).
  \end{split} \end{equation}
We turn now to the estimate of the second term in the right-hand side of \eqref{1star}. We have
\[\begin{split} 
& \sum_{\begin{subarray}{c} x,y\in \Gamma \\ y \not\sim_\epsilon x \end{subarray}}  |f_\eta(x)+f_\eta(y)|^q   p(x,y) m(y) m(x) \\
& \leq \sum_{\begin{subarray}{c} x,y\in \Gamma\\ x\sim y \\ p(x,y)m(y) \leq \epsilon \end{subarray}}|f_\eta(x)+f_\eta(y)|^q  p(x,y) m(y) m(x) + \sum_{\begin{subarray}{c}x,y\in \Gamma\\ x\sim y \\ p(x,y)m(x) \leq \epsilon \end{subarray}}|f_\eta(x)+f_\eta(y)|^q  p(x,y) m(y) m(x) \\
& \leq \epsilon \sum_{\begin{subarray}{c} x,y\in \Gamma \\ y \sim x \end{subarray} }|f_\eta(x)+f_\eta(y)|^q  [m(x)+m(y)] \\
& \leq 2^q \epsilon \left( \sum_{x\in B^\epsilon_\nu} \sum_{y\sim x} [m(x)+m(y)]  + \sum_{y\in B^\epsilon_\nu} \sum_{x\sim y} [m(x)+m(y)] \right) \\
& \quad = 2^{q+1} \epsilon \sum_{x\in B^\epsilon_\nu} \sum_{y\sim x} [m(x)+m(y)] \leq 2^{q+1} \epsilon \sum_{x\in B^\epsilon_\nu} \left[ \#\{y\in \Gamma, \, x\sim y\} m(x) + V(x,2)\right] . \\
 \end{split} \]
 Since the number of neighbors of a point is uniformly bounded by $N_{dvl}$ and since \eqref{PDV} implies that $V(x,2) \leq C_{dvl} m(x)$, one has
 \begin{equation} \label{3star} \begin{split}
 \sum_{\begin{subarray}{c} x,y\in \Gamma \\ y \not\sim_\epsilon x \end{subarray}}  |f_\eta(x)+f_\eta(y)|^q   p(x,y) m(y) m(x) & \leq 2^{q+1} \epsilon  \sum_{x\in B^\epsilon_\nu} [C_{dvl} + N_{dvl}] m(x) \\
& \leq 2^{q+1} \epsilon [C_{dvl} + N_{dvl}] m(B^\epsilon_\nu).
 \end{split}\end{equation}
Furthermore, notice that $\|f_\eta\|_{q}^q = m(B^\epsilon_\nu)$. Together with \eqref{1star}, \eqref{2star} and \eqref{3star}, we obtain
$$\|(I+P)f_\eta\|_{q}^q  \leq \left(\frac{2}{\nu} + 2^{q+1} \epsilon [C_{dvl} + N_{dvl}]\right) \|f\|_{q}^q.$$
We choose $\nu$ big enough and $\epsilon$ small enough (both depend on $\eta$ and $q$) such that $\frac{2}{\nu} + 2^{q+1} \epsilon [C_{dvl} + N_{dvl}] \leq \eta$, which gives \eqref{I+Pfeta}. The proof of $(iv) \f (vi)$ and then the proof of Theorem \ref{MainTheo} follows.
\end{proof}

\section{Application to the $L^p$-boundedness of the Riesz transform for $p\in (1,2)$ under sub-Gaussian estimates}

 \begin{defi}
 Let $(\Gamma,\mu)$ be a weighted graph and $\beta\geq 2$. We say that $(\Gamma,\mu)$ satisfies \eqref{UEm} if there exist two constants $c,C>0$ such that $p_k$ satisfies the pointwise estimates
\begin{equation} \label{UEm} \tag{UE$_\beta$}
 p_{k}(x,y) \leq \frac{C}{V(x,k^{1/\beta})} \exp\left[-c \left(\frac{d(x,y)^\beta}{k}\right)^\frac1{\beta-1}\right] \qquad  \forall x,y\in \Gamma, \, \forall k\in \N^*.
\end{equation}
\end{defi}

\begin{rmk}
The classical Gaussian estimates are ($UE_2$).
The case $\beta>2$ gives some sub-Gaussian estimates, satisfied for example, with $\beta = \log_2 5$, by the Sierpinski gasket  (see \cite{BarlowStFlour}).
\end{rmk}

Recall the definition of the positive Laplacian $\Delta :=I-P$.
We also introduce the length of the gradient 
$$\nabla f(x) : = \left( \sum_{y\in \Gamma} p(x,y) |f(x)-f(y)|^2 m(y) \right)^\frac12.$$
The operator $\nabla \Delta^{-1/2}$ is called the Riesz transform.

\medskip

With Gaussian pointwise estimates on the Markov kernel, the following result is obtained by Russ (see \cite[Theorem 1]{Russ}). The theorem in his complete form can be found in \cite[Theorem 4.2]{CCFR} (see also \cite{Fen4}).

\begin{theo} \label{TheoRieszBefore}
Let $(\Gamma, \mu)$ be a graph satisfying \eqref{DV}, \eqref{UEm} and \eqref{LB}.
For all $p\in (1,2]$, the Riesz transform $\nabla \Delta^{-1/2}$ is bounded on $L^p(\Gamma)$, that is, there exists $C_p>0$ such that 
$$\|\nabla \Delta^{-1/2} f\|_p \leq C_p \|f\|_p \qquad \forall f\in L^p(\Gamma) \cap L^2(\Gamma).$$
\end{theo}

With Theorem \ref{MainTheo}, the hypotheses in the above theorem can be weakened, and one has the following statement:

\begin{theo} \label{TheoRiesz}
Let $(\Gamma, \mu)$ be a graph satisfying \eqref{DV} and \eqref{UEm}. Assume that $P$ is analytic on $L^2(\Gamma)$.
Then, for all $p\in (1,2]$, the Riesz transform $\nabla \Delta^{-1/2}$ is bounded on $L^p(\Gamma)$, that is, there exists $C_p>0$ such that 
$$\|\nabla \Delta^{-1/2} f\|_p \leq C_p \|f\|_p \qquad \forall f\in L^p(\Gamma) \cap L^2(\Gamma).$$
\end{theo}

\begin{proof}
Let $p\in (1,2]$. Since $(\Gamma,\mu)$ satisfies \eqref{DV}, according to Theorem \ref{MainTheo}, there exists $l\in \N$, $l$ odd, such that $P^l$ satisfies \eqref{LB}. 
Since $l+1$ is even and the graph is locally doubling, $P^{l+1}$ also satisfies \eqref{LB} (see Remark \ref{P2LB}).

We define then the symmetric weight $\bar \mu$ by
$$\bar\mu_{xy} = p_{l+1}(x,y)m(x)m(y).$$
We use the notation $\bar \Gamma$ when we want to refer to $(\Gamma,\bar \mu)$. 
In particular $\bar d$, $\bar m$, $\bar \nabla$ and $\bar \Delta$ are respectively the distance, the measure, the length of the gradient, and the Laplacian in  $\bar \Gamma$.
Observe that $\bar m = m$ and $\bar d \leq d \leq (l+1)\cdot\bar d$. Therefore $\bar \Gamma$ also satisfies \eqref{DV} and \eqref{UEm}. Moreover, notice that $\bar \Delta = \Delta A$ where $A:=(I+P+\dots+P^{l})$. Theorem \ref{TheoRieszBefore} yields then 
\begin{equation} \label{BoundedRieszInter}
\left\| \bar \nabla \Delta^{-1/2} A^{-1/2} f \right\|_{p} \leq C \|f\|_{p} \qquad \forall f\in L^p(\Gamma)\cap L^2(\Gamma).
\end{equation}
Let us take a look at $A$.  First, observe that $A$ is bounded on $L^p(\Gamma)$. We want to prove that $0 \notin \Sp_{L^p}(A)$. 
Indeed, since $P^l$ satisfies \eqref{LB}, Lemma 1.3 in \cite{Dungey} implies that $\Sp_{L^p}(P^l) \subset \mathbb D \cup \{1\}$.  In particular, for any $k\in \{1,\dots,l\}$, $e^{\frac{2k\pi i}{l+1}} \notin \Sp_{L^p}(P)$, which in turn gives $0\notin \Sp_{L^p}(A)$, and so $A$ is an isomorphism on $L^p(\Gamma)$. As a consequence, the operator $A^{1/2}$ defined via functional calculus is an isomorphism on  $L^p(\Gamma)$. One can see with the spectral theorem that $A^{1/2}$ and $A^{-1/2}$ are also bounded on $L^2(\Gamma)$. Hence, \eqref{BoundedRieszInter} becomes
\begin{equation} 
\left\| \bar \nabla \Delta^{-1/2}  g \right\|_{p} \leq C \|g\|_{p} \qquad \forall g\in L^p(\Gamma)\cap L^2(\Gamma).
\end{equation}

In order to conclude the proof of Theorem \ref{TheoRiesz}, it suffices now to check that $\nabla f \leq C \bar \nabla f$ for all functions $f$ and with a constant $C$ independent of $f$, which is equivalent to proving that there exists $C>0$ such that
$p(x,y) \leq C p_{l+1}(x,y)$ for all $x,y\in \Gamma$, $x\neq y$.
Indeed, since $P^l$ satisfies \eqref{LB}, that is, since there exists $\epsilon>0$ such that $p_{l}(x,x) m(x)\geq \epsilon$ for any $x\in \Gamma$, we have
\[\begin{split}
   p(x,y) & \leq \frac1{\epsilon} p_{l}(x,x) m(x) p(x,y) \leq \frac1{\epsilon} p_{l+1}(x,y).
  \end{split}\]
  The theorem follows.
\end{proof}

\bibliographystyle{plain}
\bibliography{../Biblio}

\end{document}